\documentclass[12pt]{amsart}
\usepackage{amsmath,amsthm,amsfonts,amssymb}
\usepackage[all]{xy}

\DeclareMathOperator{\ad}{ad} \DeclareMathOperator{\Ext}{Ext}
\DeclareMathOperator{\id}{id} \DeclareMathOperator{\gp}{gp}
\DeclareMathOperator{\Hom}{Hom} \DeclareMathOperator{\im}{im}
\DeclareMathOperator{\Hopf}{Hopf} \DeclareMathOperator{\Span}{Span}

\newcommand{\coh}{{\rm H}}
\newcommand{\HH}{{\rm HH}}
\newcommand{\ra}{\rightarrow}

\newcommand{\ot}{\otimes}
\newcommand{\mtc}{\mathcal}

\numberwithin{equation}{section}
\newtheorem{lemma}[equation]{Lemma}
\newtheorem{lem}[equation]{Lemma}
\newtheorem{thm}[equation]{Theorem}

\newtheorem{cor}[equation]{Corollary}

\newtheorem{example}[equation]{Example}

\title[Hochschild cohomology]
{Hochschild cohomology of smash products and rank one Hopf algebras}
\author{S.\ M.\ Burciu}
\address{Department of Mathematics, Syracuse University,
Syracuse, NY 13210, USA} \email{smburciu@syr.edu}
\author{S.\ J.\ Witherspoon}
\address{Department of Mathematics, Texas A\&M University,
College Station, TX 77843, USA} \email{sjw@math.tamu.edu}
\thanks{The second author was partially supported by NSF grant
\#DMS-0443476 and the Alexander von Humboldt Foundation.}

\date{August 30, 2006}

\begin{document}

\maketitle

\begin{abstract}
We give some general results on the ring structure of Hochschild
cohomology of smash products of algebras with Hopf algebras. We
compute this ring structure explicitly for a large class of finite
dimensional Hopf algebras of rank one.
\end{abstract}

\section{Introduction}

Cibils and Solotar \cite{CS97} gave the ring structure of the
Hochschild cohomology of a group algebra $kG$ of a finite abelian
group $G$ over a commutative ring $k$, and Cibils \cite{C97}
conjectured a formula for this ring structure for a general finite
group $G$. Siegel and the second author \cite{SW99} proved the
conjecture. This Hochschild cohomology is a direct sum of graded
vector spaces indexed by the conjugacy classes, and cup products may
be described in terms of this decomposition via known maps and
products in group cohomology. This vector space decomposition
generalizes to a result for the Hochschild cohomology of a Hopf
algebra $H$, where the direct sum is indexed by summands of the
adjoint representation of $H$ on itself \cite[Prop.\
5.6]{GK}, but in
general there is no known  formula for the cup product in terms of
this decomposition. In the special case that $H$ is commutative, Linckelmann
\cite{Li} generalized the original
result of Cibils and Solotar.

There is a generalization in another direction, to a crossed product
of an algebra $A$ with a group algebra $kG$; again Hochschild
cohomology is a direct sum of graded vector spaces indexed by
conjugacy classes, and there is a formula for the cup product in
terms of this decomposition \cite[Thm.\ 3.16]{W04}. How much of this
theory generalizes to smash (or crossed) products with Hopf
algebras, or to Hopf Galois extensions? In this note, we begin a
program to answer this question by (1) computing the ring structure
of Hochschild cohomology for a large class of examples, namely some
finite dimensional Hopf algebras of rank one defined by Krop and
Radford \cite{KR}, by (2) giving a vector space decomposition of the
Hochschild cohomology of a smash product $A\# H$,
simultaneously generalizing the known cases $A=k$ and
$H=kG$, and by (3) giving some consequences of this decomposition in
special cases.

A rank one Hopf algebra is a generalization of a Taft algebra whose
grouplike elements may form a nonabelian group. As an algebra, it is
a smash product $B=A\# kG$, with $A=k[x]/(x^n)$ and $G$ a finite
group acting on $A$. We compute the graded vector space structure of
its Hochschild cohomology $\HH^*(B):= \Ext^*_{B\ot B^{op}}(B,B)$
(where $B^{op}$ is the algebra opposite to $B$) in Section 2 (see Theorem
\ref{gvs}). In Section 3 we use explicit chain maps first defined by
the Buenos Aires Cyclic Homology Group \cite{BuenosAires} to compute
the ring structure of $\HH^*(B)$, showing that the ring is generated
by the $G$-invariant subring $\HH^*(A)^G$ of $\HH^*(A)$ and by
$\HH^0(B)\cong Z(B)$ (see Theorem \ref{gr}).

We give our general result in Section 4 for a smash product $B=A\#
H$, where $H$ is a Hopf algebra with bijective antipode and $A$ is
an $H$-module algebra. We introduce a subalgebra $\mathcal D$ of
$B\ot B^{op}$ in (\ref{mathcald}) for which $\HH^*(B,M):=
\Ext^*_{B\ot B^{op}}(B,M)\cong\Ext^*_{\mathcal D}(A,M)$ for any
$B$-bimodule $M$ (see Theorem \ref{genthm}). If $A=k$ and $M=H$,
this yields the decomposition of $\HH^*(H)$ in terms of the adjoint
representation of $H$ on itself. If $H=kG$ and $M=B$, this yields
the decomposition indexed by conjugacy classes. In general if $M=B$,
it will give a decomposition in terms of $\mathcal D$-submodules of
$B$, and we translate the cup product on $\HH^*(B)$ to one on
$\Ext^*_{\mathcal D}(A,B)$ described explicitly at the cochain level
in (\ref{cup}). This is potentially a
first step towards understanding the cup product more directly in
terms of summands arising from the adjoint representation. It seems
difficult to generalize the next step from the special case $H=kG$,
as in this case, certain $\mathcal D$-submodules of $B$ are
coinduced from centralizer subgroups (see the proof of \cite[Lem.\
3.5]{W04}).
We also do not know if there is a more general version of our Theorem
\ref{genthm} that applies to crossed products or to Hopf Galois extensions.

In the remainder of this note, we give some consequences of Theorem \ref{genthm}.
We first return to the rank one Hopf algebras of Sections 2
and 3 and take another look at the structure of their Hochschild cohomology, this
time in relation to the adjoint representation. Next, in the special case
that $H$ is semisimple, we show that Theorem \ref{genthm} directly implies
 $\HH^*(B)\cong \HH^*(A,B)^H$, where the superscript $H$
denotes invariants (cf.\ \cite[Thm.\ 3.3]{St}), and we give some
resulting formulas for explicit cocycles (see Theorem \ref{pm}).
Finally, when $H$ is
semisimple,  another consequence of Theorem
\ref{genthm} is that the Hopf-Hochschild cohomology of $A$ introduced by
Kaygun \cite{AK} is isomorphic to the Hochschild cohomology of $B$ when
coefficients are taken in a $B$-bimodule
(Theorem \ref{Hopf-Hochschild}).
This follows from the observation that Kaygun's ``crossed product''
$A^e\rtimes H$ is isomorphic to our subalgebra $\mathcal D$ of $B^e$.

\smallskip

We work over a field $k$. For the explicit computations we require
the characteristic of $k$ to be relatively prime to the order of
$G$, however for the general results $k$ could equally well be a
commutative ring provided all algebras are projective as
$k$-modules.
Let $\ot = \ot_k$. We use
modified Sweedler notation for the coproduct $\Delta: H\rightarrow
H\ot H$ of a Hopf algebra $H$, symbolically writing
$\Delta(h)=h_1\ot h_2$ ($h\in H$).

\section{Hochschild cohomology of rank one Hopf algebras}

Let $G$ be a finite group whose order is relatively prime to the
characteristic of $k$. Let $\chi: G\rightarrow k^{\times}$ be a
character, that is a group homomorphism from $G$ to the
multiplicative group of $k$. Let $n\geq 2$ be a positive integer and
$A=k[x]/(x^n)$. Then $G$ acts by automorphisms on $A$ via
$$
   {}^g x = \chi(g) x
$$
for all $g\in G$. Let $B=A\# kG$, the corresponding {\em skew group
algebra} (or {\em smash product} of $A$ and $kG$): As a vector
space, $B=A\ot kG$, and the multiplication is
$$
  (a\ot g)(b\ot h) = a ({}^g b)\ot gh
$$
for all $a,b\in A$ and $g,h\in G$. We abbreviate $a\ot g$ by $ag$.

Assume there is a central element $g_1\in G$ such that $\chi(g_1)$
is a primitive $n$th root of 1. Then $B$ is a Hopf algebra with
coproduct $\Delta$ defined by
\begin{equation}\label{coproducts}
  \Delta(x)= x\otimes 1 + g_1\otimes x \ \ \ \mbox{ and } \ \ \ \Delta(g)
   =g\otimes g,
\end{equation}
counit $\varepsilon$ by $\varepsilon(x)=0$ and $\varepsilon(g)=1$,
and antipode $S$ by $S(x)=-g_1^{-1}x$ and $S(g)=g^{-1}$, for all
$g\in G$. This generalization of a Taft algebra is defined in
\cite{AS} for abelian groups $G$, and generalized further in
\cite{KR} (but with the opposite coproduct). Note that we do not use
the coalgebra structure of $B$ until Section 4.

In order to compute the Hochschild cohomology of $B$, we use the
following subalgebra of $B^e=B\otimes B^{op}$:
\begin{equation}\label{mathcald-gp}
  {\mathcal D} := A^e\# kG\cong \bigoplus_{g\in G} (Ag\ot Ag^{-1})
   \subset B^e,
\end{equation}
where the action of $G$ on $A^e$ is diagonal, that is ${}^g(a\ot b)
={}^ga\ot {}^gb$. The indicated isomorphism is given by $(a\otimes
b)g \mapsto ag\otimes ({}^{g^{-1}}b)g^{-1}$ for all $a,b\in A$ and
$g\in G$. Note that $A$ is a $\mathcal D$-module under left and
right multiplication.
The algebra $\mathcal D$ is sometimes denoted $\Delta$ in the literature
on group-graded algebras.

It is known that the Hochschild cohomology
$\HH^*(B):=\Ext^*_{B^e}(B,B)$ satisfies
\begin{equation}\label{ESL}
  \HH^*(B) \cong \Ext^*_{\mathcal D}(A,B)
\end{equation}
as graded algebras. This is a consequence of the Eckmann-Shapiro
Lemma and the isomorphism of $B^e$-modules,
$\displaystyle{B\stackrel{\sim}{\rightarrow} A\uparrow_{\mathcal
D}^{B^e} :=B^e\ot_{\mathcal D}A}$ given by $b\mapsto (b\ot 1)\ot 1$
with inverse $(b\ot c)\ot a\mapsto bac$. (See for example
\cite[Lemma 3.5]{W04}, valid more generally for some crossed
products.) Alternatively (\ref{ESL}) follows from our generalization
to smash products with Hopf algebras, Theorem \ref{genthm} below. As
the characteristic of $k$ is relatively prime to $|G|$, there is a
further isomorphism $\Ext^*_{\mathcal D}(A,B)\cong
\Ext^*_{A^e}(A,B)^G=\HH^*(A,B)^G$ where the latter consists of
invariants under the action induced from that of $G\subset {\mathcal
D}$ on ${\mathcal D}$-modules. (The resulting isomorphism
$\HH^*(B)\cong \Ext^*_{A^e}(A,B)^G$ also follows from \cite[Cor.\
3.4]{St} or from (\ref{ss}) below.) Again as the characteristic of
$k$ is relatively prime to $|G|$, $G$-invariants may be taken in a
complex prior to taking cohomology. This will be our approach in
proving the following theorem for a rank one Hopf algebra $B=A\#
kG$.

\begin{thm} \label{gvs}
Let $N=\ker\chi\subset G$. For all $i\geq 0$,
$$
  \HH^{2i}(B)\cong \HH^{2i+1}(B)\cong
    \left\{\begin{array}{cl} Z(kN), & \mbox{ if }\chi^{in}=1\\
                                         0,   & \mbox{ otherwise}.
                      \end{array}\right.
$$
\end{thm}

\begin{proof}
The following is an $A^e$-free resolution of $A$ \cite[Exer.\
9.1.4]{weibel}:
\begin{equation}\label{projres}
\cdots \stackrel{\cdot v}{\longrightarrow} A^e
       \stackrel{\cdot u}{\longrightarrow} A^e
       \stackrel{\cdot v}{\longrightarrow} A^e
       \stackrel{\cdot u}{\longrightarrow} A^e
       \stackrel{\mu}{\longrightarrow} A \rightarrow 0,
\end{equation}
where $u=x\ot 1-1\ot x$, $v=x^{n-1}\ot 1 + x^{n-2}\ot x +\cdots +
1\ot x^{n-1}$, and $\mu$ is multiplication.
This becomes a $\mathcal D$-projective resolution of $A$ as follows.
The action of $G\subset {\mathcal D}$ in degree 0 is diagonal on
$A^e$, $g\cdot (a\ot b) = {}^g a\ot {}^g b$ for all $a,b\in A$ and
$g\in G$. In all other degrees the action must be modified in order
that the maps $\cdot u$ and $\cdot v$ be maps of $\mathcal
D$-modules. In degree $2i$, $g\cdot (a\ot b) =\chi(g)^{in}
({}^ga)\ot {}^gb$, and in degree $2i+1$, $g\cdot (a\ot b) =
\chi(g)^{in+1} ({}^g a)\ot {}^gb$, for all $a,b\in A$ and  $g\in G$. With
these actions, (\ref{projres}) is indeed a $\mathcal D$-projective
resolution of $A$: Since the characteristic of $k$ does not divide
the order of $G$, a $\mathcal D$-module is projective if and only if
its restriction to $A^e$ is projective. (An $A^e$-splitting map of
$\mathcal D$-modules may be ``averaged'' by applying
$\frac{1}{|G|}\sum _{g\in G} g$ to obtain a $\mathcal D$-splitting
map.)

According to the isomorphism (\ref{ESL}), we must now apply
$\Hom_{\mathcal D}( - , B)$ to (\ref{projres}). Now  $\Hom_{\mathcal
D}(A^e,B)\cong \Hom_{A^e}(A^e,B)^G$ where $G$ acts on
$\Hom_{A^e}(A^e,B)$ by $(g\cdot f)(a\otimes b) =g f(g^{-1} \cdot
(a\otimes b))g^{-1}$ for all $f\in\Hom_{A^e}(A^e,B), g\in G$. Such a
homomorphism is determined by its value on $1\ot 1$. We identify $B$
with $\Hom_k(k,B)\cong \Hom_{A^e}(A^e,B)$, under the correspondence
$b\mapsto f_b$ where $f_b(1\ot 1)=b$. Thus applying $\Hom_{\mathcal
D}(-,B)$ to (\ref{projres}) yields the complex
\begin{equation}\label{bg}
  \cdots \stackrel{(\cdot v)^*}{\longleftarrow\joinrel\relbar} B^G
         \stackrel{(\cdot u)^*}{\longleftarrow\joinrel\relbar} B^G
         \stackrel{(\cdot v)^*}{\longleftarrow\joinrel\relbar} B^G
         \stackrel{(\cdot u)^*}{\longleftarrow\joinrel\relbar} B^G
           \leftarrow 0,
\end{equation}
the action of $G$ on $B$ depending on the degree as stated above. In
degree $2i$, this action is
$$
  (g\cdot f_b)(1\ot 1) = gf_b(g^{-1}\cdot (1\ot 1))g^{-1} =
    g f_b(\chi(g)^{-in}\ot 1) g^{-1} = \chi(g)^{-in} gbg^{-1},
$$
so that $g\cdot b = \chi(g)^{-in} gbg^{-1}$ for all $g\in G$ and $b\in
B$. Similarly, in degree $2i+1$, $g\cdot b =\chi(g)^{-in-1}
gbg^{-1}$. Thus in degree $2i$, as
$\chi(g_1)$ is a primitive $n$th root of 1, if $x^j$ is invariant
then $x^j=g_1\cdot x^j =\chi(g_1)^{-in+j}x^j=\chi(g_1)^jx^j$,
implying $j=0$. It follows that $B^G\subset kG$. Applying the
formula for the action in degree $2i$, we find $B^G=Z(kG)$ (the
center of the group algebra $kG$) in case $\chi^{in}=1$, and
otherwise there are no invariants. Thus in degree $2i$, $B^G$ is
spanned by elements of the form $\sum_{h\in G/C(g)} hgh^{-1}$ in
case $\chi^{in}=1$, where $C(g)$ is the centralizer of an element $g$ in $G$.
Similarly, in degree $2i+1$, the invariants are spanned by elements
of the form $\sum_{h\in G/C(g)} xhgh^{-1}$ in case $\chi^{in}=1$,
and otherwise there are no invariants.

The maps $(\cdot u)^*$ and $(\cdot v)^*$ are:
\begin{eqnarray*}
  (\cdot u)^*(b) \! &=& \! (\cdot u)^* f_b(1\ot 1) = f_b(u)  = xb-bx,\\
  (\cdot v)^*(b) \!&=&\! x^{n-1}b +x^{n-2}bx +\cdots + bx^{n-1},
\end{eqnarray*}
for all $b\in B^G$. In particular $(\cdot v)^*$ is the $0$-map:
$(\cdot v)^* (\sum_{h\in G/C(g)} xhgh^{-1}) = 0$ since $x^n=0$. Thus
$\ker (\cdot v)^* =B^G$ and $\im (\cdot v)^* =0$. In degree $2i$, if
$B^G\neq 0$, then $\ker(\cdot u)^*$ is spanned by those conjugacy
class sums $\sum_{h\in G/C(g)}hgh^{-1}$ for which $g\in\ker \chi =
N$. In degree $2i+1$, $\im (\cdot u)^*$ is spanned by those elements
$\sum_{h\in G/C(g)} xhgh^{-1}$ for which $g\not\in N$. Therefore
\begin{eqnarray*}
  \HH^{2i}(B)\! &= &\! \ker (\cdot u)^* /\im (\cdot v)^* \cong Z(kN)
   \mbox{ if }\chi^{in}=1, \mbox{ and 0 otherwise},\\
\HH^{2i+1}(B) \! & = &\! \ker (\cdot v)^* /\im (\cdot u)^* \cong
xZ(kN) \mbox{ if }
   \chi^{in}=1, \mbox{ and 0 otherwise}.
\end{eqnarray*}
As vector spaces, we thus have $\HH^{2i}(B)\cong \HH^{2i+1}(B)$ for
each $i$.
\end{proof}

\section{The ring structure}

We next compute the ring structure of $\HH^*(B)$, where $B=A\# kG$
is the rank one Hopf algebra defined in Section 2. In order to do
so, we compare the resolution (\ref{projres}) with the bar
resolution of $A$,
\begin{equation}\label{barres}
  \cdots\stackrel{\delta_3}{\longrightarrow}
  A^{\ot 4}\stackrel{\delta_2}{\longrightarrow}
  A^{\ot 3}\stackrel{\delta_1}{\longrightarrow}
  A^{e} \stackrel{\mu}{\longrightarrow} A\rightarrow 0,
\end{equation}
where $\delta_i(a_0\ot a_1\ot\cdots\ot a_{i+1}) =\sum_{j=0}^i (-1)^j
    a_0\ot\cdots\ot a_ja_{j+1}\ot\cdots\ot a_{i+1}$
for $a_0,\ldots,a_{i+1}\in A$.
The cup product on $\HH^*(B)$ is defined at the chain level, with respect
to the bar resolution (\ref{barres}) with $A$ replaced by $B$, as
follows: Let $f\in \Hom_{B^e}(B^{\ot (l+2)},B)\cong \Hom_k(B^{\ot
l},B)$ and $f'\in\Hom_{B^e}(B^{\ot (m+2)},B)\cong \Hom_k(B^{\ot
m},B)$. Then
\begin{equation}\label{cupprodformula}
    (f\smile f')(b_1\ot \cdots\ot b_{l+m})=f(b_1\ot\cdots\ot b_l)
    f'(b_{l+1}\ot\cdots\ot b_m)
\end{equation}
for all $b_1,\ldots,b_{l+m}\in B$. It is convenient to consider an
$l$-cochain sometimes to be an element of $\Hom_{B^e}(B^{\ot
(l+2)},B)$ and other times to be an element of $\Hom_k(B^{\ot
l},B)$. This should cause no confusion.

Consider the cup product on $\HH^*(A,B)$ given by (\ref{cupprodformula}),
where $B^{\ot l}$, $B^{\ot m}$ are replaced by $A^{\ot l}$, $A^{\ot m}$.
The isomorphism $\HH^*(B)\cong \HH^*(A,B)^G$ described in Section 2
preserves cup products,
that is $\HH^*(B)$ is isomorphic to the $G$-invariant subalgebra of
$\HH^*(A,B)$. This is known, and also follows from our more
general results at the end of Section 4 as $kG$ is semisimple, but
we outline a direct proof using the algebra $\mathcal D$ in this
case. Note that the bar resolution for $B$ (as $B^e$-module) is
induced from the $\mathcal D$-projective resolution of $A$:
\begin{equation}\label{dres}
\cdots \stackrel{\delta_3}{\longrightarrow} {\mathcal D}_2
 \stackrel{\delta_2}{\longrightarrow} {\mathcal D}_1
 \stackrel{\delta_1}{\longrightarrow} {\mathcal D}_0
 \stackrel{\mu}{\longrightarrow} A \rightarrow 0
\end{equation}
where ${\mathcal D}_0 = {\mathcal D}$ and
$${\mathcal D}_m =\Span_k \{ a_0g_0\ot \cdots \ot a_{m+1}g_{m+1} \mid
 a_i\in A , g_i\in G, g_0\cdots g_{m+1} = 1\}
$$
is a $\mathcal D$-submodule of $B^{\ot (m+2)}$. An isomorphism
${\mathcal D}_m\uparrow_{\mathcal D}^{B^e}:=
\displaystyle{B^e\ot_{\mathcal D} {\mathcal D}_m \stackrel
{\sim}{\rightarrow} B^{\ot (m+2)}}$ is given by
$$
  (b_{-1}\ot b_{m+2})\ot (b_0\ot\cdots \ot b_{m+1})
\mapsto b_{-1}b_0\ot b_1\ot \cdots\ot b_{m+1}b_{m+2},
$$
and its inverse by
$$
  a_0g_0\ot a_1g_1\ot \cdots \ot a_{m+1}g_{m+1}\mapsto
  (1\ot g_0\cdots g_{m+1})\ot (a_0g_0\ot\cdots\ot a_{m+1} g_m^{-1}\cdots
   g_0^{-1}).$$
The bar resolution (\ref{barres}) for $A$ is compatible with the
action of ${\mathcal D}=A^e\# kG$ given by the usual action of $A^e$
and the diagonal action of $G$ on tensor products $A^{\ot m}$. Thus
(\ref{barres}) is in fact a $\mathcal D$-projective resolution of
$A$. 
There is a $\mathcal D$-map from (\ref{dres}) to (\ref{barres}) given by
$$
  a_0g_0\ot\cdots\ot a_{m+1}g_{m+1} \mapsto
    a_0\ot {}^{g_0}a_1\ot {}^{g_0g_1}a_2\ot\cdots\ot {}^{g_0\cdots g_m}a_{m+1}
$$
for all $a_i\in A$ and $g_i\in G$ (see \cite[(5.2)]{CGW}).
Under this map and the identification $B^{\ot (m+2)}\cong 
{\mathcal D}_m\uparrow_{\mathcal D}^{B^e}$, it can be seen that the
cup product (\ref{cupprodformula}) on $\HH^*(B)$ indeed corresponds to that on
$\HH^*(A,B)^G\subset \HH^*(A,B)$ induced by multiplication on $B$.

We first need a chain map $\phi_*$ from (\ref{projres}) to the bar
complex (\ref{barres}) for $A$. This was found in a more general
setting in \cite{BuenosAires}. We give the maps explicitly in our
setting.
Define $\phi_m : A^e\rightarrow A^{\ot (m+2)}$ by
$$\phi_{2l}(1 \ot 1)=1 \ot \alpha_l \ \ \mbox{ and } \ \ \phi_{2l+1}(1 \ot
1)=1\ot x \ot \alpha_l
$$
where $\alpha_0=1$ and if $l\geq 1$,
$$\alpha_l=\sum_{\{i_1+i_2+ \cdots +i_{l+1}=ln-l | i_1, i_2, \cdots , i_l \geq 1\}}x^{i_1}\ot x \ot x^{i_2} \ot x\ot \cdots \ot x \ot x^{i_{l+1}}.$$
In this formula, we emphasize that $i_{l+1}$ is allowed to be $0$
whereas $i_1,i_2,\ldots,i_{l}$ must be greater than 0. Note that our
maps differ from those in \cite{BuenosAires} by a sign due to
our choice $u=x\ot 1 - 1\ot x$.

We next define a chain map $\psi_*$ from (\ref{barres}) to
(\ref{projres}) as in \cite{BuenosAires}.
Define $\psi_{2l}: A^{\ot (2l+2)} \ra A^e$ and $\psi_{2l+1}:A^{\ot
(2l+3)}\rightarrow A^e$ by
\begin{eqnarray*}
\psi_{2l}(1\! \ot\! x^{i_1}\! \ot\! x ^{i_2}\!\ot\!  \cdots \!\ot\!
x^{i_{2l}}\!\ot\! 1)\!\! &=&\!\! 1\!\ot\!
x^{i_1+i_2-n}x^{i_3+i_4-n}\! \cdots\!
x^{i_{2l-1}+i_{2l}-n},\\
\psi_{2l+1}(1 \!\ot\! x^{i_1}\!\ot\! x ^{i_2}\!\ot\!  \cdots \!\ot\!
x^{i_{2l+1}}\!\ot\! 1)\!\! &=&\!\!\sum_{m=0}^{i_1-1} x^m\!\ot\!
x^{i_1-m-1} x^{i_2+i_3-n}x^{i_4+i_5-n}\! \cdots\!
x^{i_{2l}+i_{2l-1}-n},
\end{eqnarray*}
where $x^j$ is defined to be 0 if $j<0$. By \cite[Prop.\
1.5]{BuenosAires}, $\phi_*$ and $\psi_*$ are indeed chain maps.
Further, both $\phi_*$ and  $\psi_*$ are compatible with the action
of $G$ and thus they are $\mathcal D$-maps.

\begin{thm} \label{gr}
Let $N=\ker\chi \subset G$ and let $p$ be the order of $\chi^n$ in
$\Hom_{\gp}(G,k^{\times})$. Then there is an isomorphism of graded
algebras
$$
  \HH^*(B)\cong Z(kN)\ot k[y,z]/(z^2),
$$
where $\deg y =  2p$ and $\deg z = 1$. In particular $\HH^*(B)$ is
generated by $\HH^*(A)^G\cong k[y,z]/(z^2)$ and $\HH^0(B)\cong
Z(kN)$.
\end{thm}

\begin{proof}
Computations will be done at the chain level using the complex
(\ref{bg}). First let $a, b \in B^G$ be two elements of degrees $2l$
and $2m$, respectively, and let $f_a$ and $f_b$ be the corresponding
functions from $A^e$ to $B^G$. Applying the chain maps $\phi_*$,
$\psi_*$ given above and the definition of cup product on the bar
resolution, the cup product $f_a \smile f_b$ is defined by
$$(\psi_{2l}^*f_a \smile \psi_{2m}^*f_b)\phi_{2l+2m}(1 \ot
1)\hspace{3in}$$
$$\hspace{1in}=(\psi_{2l}^*f_a \smile \psi_{2m}^*f_b)(\sum 1\ot x^{i_1}\ot x \ot
x^{i_2} \ot x\ot \cdots \ot x \ot x^{i_{l+m+1}})$$ where the sum is
over all indices $i_s$ such that $i_1 \cdots +i_{l+m+1}=(l+m)(n-1 )$
and $ i_1, i_2, \cdots , i_{l+m} \geq 1$. Identifying $f_a,f_b$
with $a,b$, we have
\begin{eqnarray*}
a \smile b&=&\sum \psi_{2l}^*f_a (x^{i_1}\ot x\ot \cdots \ot
x^{i_{l}}\ot x) \psi_{2m}^*f_b(x^{i_{l+1}}\ot x \ot\cdots \ot
x^{i_{l+m}}\ot x)x^{i_{l+m+1}}\\
&=&\sum f_a(1\ot x^{i_1+1-n}\cdots x^{i_l+1-n})f_b(1\ot
x^{i_{l+1}+1-n}\cdots
  x^{i_{l+m}+1-n}) x^{i_{l+m+1}} \
= \ ab
\end{eqnarray*}
since the sum has only one nonzero term, the one where
$i_1=i_2=\cdots=i_{l+m}=n-1$ and $i_{l+m+1}=0$.
If $a,b \in B^G$ are elements of degrees $2l$ and $2m+1$,
respectively, then a similar calculation shows that $a\smile b = ab$.

Finally, let $a,b \in B^G$ be  elements of degrees $2l+1$
and $2m+1$, respectively. Then the cup product $f_a \smile f_b$ is
given by
$$(\psi_{2l+1}^*f_a \smile \psi_{2m+1}^*f_b)\phi_{2l+2m+2}(1 \ot
1)\hspace{2.5in}$$
$$\hspace{1in}=(\psi_{2l+1}^*f_a \smile \psi_{2m+1}^*f_b)(\sum  1\ot x^{i_1}\ot
x\ot
 \cdots \ot x \ot x^{i_{l+m+2}})$$ where the sum is
over all indices $i_s$ such that $i_1 +\cdots
+i_{l+m+2}=(l+m+1)(n-1)$ and $ i_1, i_2, \cdots ,i_{l+m+1} \geq
1$. So $a\smile b$ is
\begin{eqnarray*}
\!\! &=&\!\! \sum \psi_{2l+1}^*f_a ( x^{i_1}\!\ot\! x  \!\ot\!\cdots\! \ot\!
x^{i_{l+1}})\psi_{2m+1}^*f_b(x\!\ot
\! x^{i_{l+2}}\! \ot\! \cdots\! \ot\! x^{i_{l+m+1}}\!\ot\! x)x^{i_{l+m+2}} \\
\!\! &=&\!\! \sum \sum_{j=0}^{i_1-1} f_a(x^j\!\ot\! x^{i_1-j-1}x^{i_2+1-n}\!\cdots\!
   x^{i_{l+1}+1-n}) f_b(1\!\ot\! x^{i_{l+2}+1-n}\!\cdots\!
   x^{i_{l+m+1}+1-n})x^{i_{l+m+2}}\\
 \!\!&=&\!\! \sum_{i_1\geq 1 , \ i_1+i_{l+m+2}=n-1} \ \sum_{j=0}^{i_1-1}
    f_a(x^j\ot x^{i_1-j-1}) f_b(1\ot 1) x^{i_{l+m+2}}\\
  \!\!&=&\!\! \sum \sum x^j a x^{i_1-j-1} b x^{i_{l+m+2}} \ = \ 0
\end{eqnarray*}
since $a,b\in kGx$, $i_1 + i_{l+m+2}=n-1$, and $x^n=0$.

Now let $z=x$ in degree 1 and $y=1$ in degree $2p$. Comparing to the
proof of Theorem \ref{gvs}, we see that $y$ and $z$ together with
$\HH^0(B)\cong Z(kN)$ generate $\HH^*(B)$, and the ring structure is
as claimed.
\end{proof}

Note that $\HH^*(B)$ is finitely generated. It has been conjectured
that the Hochschild cohomology ring, modulo nilpotent elements, of
any finite dimensional algebra is finitely generated \cite{SS}.

We remark that the cup products could equally well have been computed
using Yoneda composition. The resolution (\ref{projres}) is a
$\mathcal D$-projective resolution of $A$, and may be induced to
$B^e$ to obtain a $B^e$-projective resolution of $A\uparrow
_{\mathcal D}^{B^e}\cong B$. The technique for computing Yoneda
compositions from a projective resolution given in
\cite[\S2.6]{benson91a} applies to this resolution to yield an
alternative proof of Theorem \ref{gr}.

It would be interesting to determine the Hochschild cohomology more
generally for {\em all} finite dimensional rank one Hopf algebras,
including those for which the relation $x^n=0$ is replaced by
$x^n=g_1^n-1$ (see \cite{KR}). This would require a different
approach. A generalization in another direction would be to allow
the characteristic of $k$ to divide $|G|$ while remaining relatively
prime to $n$.

\section{Hochschild cohomology of smash products}

In this section, we let $A$ be any $k$-algebra and $H$ any Hopf
algebra over $k$, with bijective antipode $S$, for which $A$ is an
$H$-module algebra. That is, $A$ is an $H$-module for which ${}^h
(ab) = ({}^{h_1} a) ({}^{h_2} b)$ and ${}^h 1=\varepsilon(h)$ for
all $a,b\in A$ and $h\in H$. Let $B=A\# H$ be the {\em smash product} of
$A$ with $H$: As a vector space, $B=A\ot H$, and multiplication is
given by
$$
   (a\ot h)(b\ot l)=a( {}^{h_1} b) \ot h_2 l
$$
for all $a,b\in A$ and $h,l\in H$. We abbreviate $a\ot h$ by $ah$.

We now generalize the algebra ${\mathcal D}$ defined in (\ref{mathcald-gp}) in
the case $H=kG$.
Let $\delta: H\rightarrow H\ot H^{op}$ be the map given by
$\delta(h)= h_1\ot S(h_2)$ for all $h\in H$.
Note that $\delta$ is injective as its composition with $\id\ot\varepsilon$
is injective, so that $H\cong \delta(H)$.
Let
\begin{equation}\label{mathcald}
{\mathcal D}:=(A\ot A^{op}) \delta(H),
\end{equation}
 a subalgebra of $B^e$:
To see that $\mathcal D$ is closed under multiplication, use the
relation $aS(h) = S(h_1)({}^{h_2}a)$ for all $h\in H, a\in A$:
\begin{eqnarray*}
  (a\ot b) (h_1\ot S(h_2))(c\ot d)(l_1\ot S(l_2)) &=&
     (a\ot b)(h_1 c\ot d S(h_2)) (l_1\ot S(l_2))\\
   &=& (a\ot b) (({}^{h_1} c)h_2\ot S(h_3) ({}^{h_4}d))(l_1\ot S(l_2))\\
   &=& (a\ot b)({}^{h_1}c\ot {}^{h_4}d)(h_2\ot S(h_3))(l_1\ot S(l_2))
\end{eqnarray*}
for all $a,b,c,d\in A$ and $h,l\in H$.
Unlike the case $H=kG$, the algebra $\mathcal D$ appears not to be a smash 
product in general.

Note that $A$ is a $\mathcal D$-module under left and right
multiplication since $h_1 a S(h_2)=({}^{h_1} a) h_2 S(h_3) = {}^{h}
a\in A$ for all $h\in H$ and $a\in A$. Let $A\uparrow_{\mathcal
D}^{B^e}=B^e\ot_{\mathcal D}A$ denote the induced left $B^e$-module.

\begin{lemma}\label{lem2}
There is an isomorphism of left $B^e$-modules, $B\cong
A\uparrow_{\mathcal D}^{B^e}$.
\end{lemma}

\begin{proof}
First note that $H^e= (H\ot 1)\delta(H)$ as sets: If $a,b\in H$, we
have
$$
  a\ot S(b) = a\varepsilon(b_1)\ot S(b_2)
   = aS(b_1)b_2\ot S(b_3)
   = (aS(b_1)\ot 1) (b_2\ot S(b_3)),
$$
an element of $(H\ot 1) \delta(H)$. This suffices since $S$ is
bijective.

Now define $B^e$-maps $\phi$ and $\psi$:
$$
\begin{array}{ccccccc}
  B^e\ot_{\mathcal D} A &\stackrel{\phi}{\relbar\joinrel\longrightarrow}
  & B & \ \ \ ,\hspace{.5in}
   & B &\stackrel{\psi}{\relbar\joinrel\longrightarrow}
  & B^e\ot_{\mathcal D} A   \ \ \ .\\
 (b\ot c)\ot a & \mapsto & bac & & b &\mapsto & (b\ot 1)\ot 1
\end{array}
$$
That $\psi$ is a $B^e$-map uses the relation from the first
paragraph. Clearly $\phi$ is well-defined. We next check that $\phi$
and $\psi$ are inverses. By the above arguments, $
  B^e = H^e A^e = (H\ot 1) \delta(H) A^e.
$ Note that ${\mathcal D} = A^e\delta(H) = \delta(H) A^e$: If $h\in
H$ and $a,b\in A$, then $$(h_1\ot S(h_2))(a\ot b) = h_1a\ot bS(h_2) =
(h_1.a)h_2\ot S(h_3)h_4.b = (h_1.a\ot h_4.b)(h_2\ot S(h_3)),$$ so
that $\delta(H)A^e\subseteq A^e\delta(H)$. The other containment may
be shown similarly. 

We claim that $B^e$ is a free right $\mathcal D$-module, with free
$\mathcal D$-basis given by any $k$-basis of $H\ot 1$. In case
$A=k$, this follows from the fact that a tensor product of a free
$H$-module with another  module is free \cite[Prop.\
3.1.5]{benson91a}. In general, note that $B^e$ is free over $A^e$
with free basis any $k$-basis of $H^e$, by construction, and we may
take as $k$-basis of $H^e$ the product of a $k$-basis of $H\ot 1$
with a $k$-basis of $\delta(H)$, again by \cite[Prop.\ 3.1.5]{benson91a}. 
This shows that a $k$-basis of
$H\ot 1$ forms a free $\mathcal D$-basis of $B^e$. Therefore we may
write elements of $B^e\ot_{\mathcal D} A$ as linear combinations of
elements $(h\ot 1)\ot a$ ($h\in H, a\in A$). Then
$$
  (\psi\circ\phi)((h\ot 1)\ot a) = \psi(ha)=(ha\ot 1)\ot 1 =(h\ot 1)\ot a
$$
for all $h\in H, a\in A$, since $a\ot 1\in {\mathcal D}$, and $
  (\phi\circ\psi)(b) =\phi((b\ot 1)\ot 1) = b
$ for all $b\in B$.
\end{proof}

The following theorem generalizes part of \cite[Lemma 3.5]{W04}.

\begin{thm}\label{genthm}
Let $M$ be a $B$-bimodule, and $\mathcal D$ the subalgebra of $B^e$
defined in (\ref{mathcald}). Then
$$\HH^*(B,M)\cong \Ext^*_{\mathcal D}(A,M)$$
as graded vector spaces.
\end{thm}

We remark that any decomposition of $M$ into a direct sum of
$\mathcal D$-submodules now leads to a similar decomposition of
$\HH^*(B,M)$.

\begin{proof}
Since $B^e$ is a free right $\mathcal D$-module we may apply the
Eckmann-Shapiro Lemma  and Lemma \ref{lem2} to obtain the claimed
isomorphism.
\end{proof}

\begin{cor} \label{cor}
The graded vector space $\Ext^*_{\mathcal D}(A,B)$ has a ring
structure for which
$$\HH^*(B) \cong \Ext^*_{\mathcal D}(A,B)$$ as graded algebras.
\end{cor}

As a consequence, $\HH^*(B)$ has a graded vector space decomposition
indexed by $\mathcal D$-summands of $B$. In case $B=kG$, there are
$\mathcal D$-summands that are indexed by the conjugacy classes of
the group, and this leads to a useful description of cup products on
$\HH^*(kG)$ that was used to compute several examples \cite{SW99}.
In case $B=H$ is a commutative Hopf algebra, $H$ is trivial as a
$\mathcal D$-module. Combined with an explicit description of the
cup product given in (\ref{cup}) below, this yields an isomorphism of graded
algebras $\HH^*(H)\cong H\ot \Ext^*_H(k,k)$, an alternative proof of
\cite[Thm.\ 1]{Li}. In case $B=A\# H$ with $H=kG$, the cup product
on $\Ext^*_{\mathcal D}(A,B)$ is described in \cite[Thm.\ 3.16]{W04}
in terms of summands indexed by conjugacy classes.

We give an explicit formula for the cup product on $\Ext^*_{\mathcal
D} (A,B)$ referred to in Corollary \ref{cor}, by expressing the bar
resolution of $B$ as induced from a $\mathcal D$-resolution of $A$.
Define $\delta^m:H^{\ot (m+1)}\rightarrow H^{\ot (m+2)}$ by
$$
  \delta^m(h^{0}\ot\cdots\ot h^{m})=h_1^{0}\ot\cdots\ot h_1^{m}
  \ot S(h_2^{0}\cdots h_2^{m})
$$
for $h^{i}\in H$. Let
$$
  {\mathcal D}_m:=(A^{\ot (m+2)})\delta^m(H^{\ot (m+1)}),
$$
where indicated products occur in $B$.
Note that $\delta^0=\delta$ and ${\mathcal D}_0={\mathcal D}$
since $$ah_1\ot bS(h_2) = ah_1\ot S(h_2)({}^{h_3}b)=(a\ot {}^{h_3}b)
(h_1\ot S(h_2))\in {\mathcal D}\subset B^e.$$
By Lemma \ref{lem2} we have an isomorphism of $B^e$-modules,
$A\uparrow_{\mathcal D}^{B^e}\cong B$. By construction, ${\mathcal
D}_0\uparrow_{\mathcal D}^{B^e}\cong B^e$. A calculation shows that
for all $m\geq 0$, $\displaystyle{{\mathcal D}_m\uparrow_{\mathcal
D}^{B^e} \stackrel{\sim}{\rightarrow}B^{\ot (m+2)}}$ as
$B^e$-modules, via the map
\begin{equation}\label{map}
   (b_{-1}\ot b_{m+2})\ot (b_0\ot\cdots
\ot  b_{m+1}) \mapsto b_{-1}b_0\ot b_1\ot\cdots\ot b_m\ot
  b_{m+1}   b_{m+2},
\end{equation}
whose inverse is
\begin{equation}\label{inverse}
  a_0h^{0}\! \ot\! a_1h^{1}\!\ot\!\cdots\!\ot\! a_{m+1}h^{m+1}\mapsto
  (1\!\ot\! h_3^{0}\!\cdots\! h_3^{m}h^{m+1})\!\ot\!
  (a_0h_1^{0}\!\ot\!\cdots\!\ot\! a_m h_1^{m}\!\ot\! a_{m+1}
   S(h_2^{0}\!\cdots\! h_2^{m})).
\end{equation}
We claim that ${\mathcal D}_m$ is $\mathcal D$-projective: First
note that $${\mathcal D}_m ={\mathcal D}\{1\ot a_1h_1^1\ot\cdots\ot
a_mh_1^m\ot S(h_2^1\cdots h_2^m)\}$$ as sets.
We use this to define an isomorphism $\displaystyle{\mathcal
D}_m\stackrel{\sim}{\rightarrow} (A\ot B^{\ot m}\ot A)\uparrow
^{\mathcal D}_{A^e}$ via the map 
$$(a_0l_1\ot a_{m+1}S(l_2))(1\ot a_1h_1^1\ot\cdots\ot a_mh_1^m\ot S(h_2^1\cdots h_2^m))
\hspace{1in}$$
$$\hspace{1in}\mapsto (a_0l_1\ot a_{m+1}S(l_2))\ot (1\ot a_1h^1\ot\cdots\ot
  a_mh^m\ot 1)$$
whose inverse is
$$
  (a_{-1}l_1\ot a_{m+2}S(l_2))\ot a_0\ot a_1h^1\ot\cdots\ot a_mh^m\ot a_{m+1}
\hspace{1in}$$
$$\hspace{1in}\mapsto a_{-1}l_1a_0\ot a_1h_1^1\ot\cdots\ot a_m h_1^m\ot S(h_2^1\cdots h_2^m)
  a_{m+1}a_{m+2}S(l_2).
$$
Since $A\ot B^{\ot m}\ot A$ is clearly $A^e$-free and $\mathcal D$
is free as a right $A$-module, the induced $\mathcal D$-module ${\mathcal D}_m$ is
$\mathcal D$-free. The differentials for the bar resolution of
$B$ preserve ${\mathcal D}_*$, and it may be checked that ${\mathcal
D}_*$ is a resolution under the restriction of these differentials.
The bar resolution of $B$ is thus induced from ${\mathcal D}_*$. If
$f:{\mathcal D}_l\rightarrow B$ and $f':{\mathcal D}_m\rightarrow B$
are two cocycles, define $f\smile f': {\mathcal D}_{l+m}\rightarrow
B$ by
\begin{equation}\label{cup}
  (f\!\smile\! f')(a_0h_1^{0}\!\ot\!\cdots\!\ot\! a_{l+m}h_1^{l+m}\!\ot\! a_{l+m+1}
    S(h_2^{0}\!\cdots\! h_2^{l+m}))\!
   =\! f(a_0h_1^{0}\!\ot\!\cdots\!\ot\! a_lh_1^{l}\!\ot\! S(h_2^{0}\!\cdots\! h_2^{l}))
   \cdot
\end{equation}
$$
  \hspace{1cm}  f'( h_3^{0}\cdots h_3^{l}
\ot a_{l+1}h_1^{l+1}\ot\cdots\ot a_{l+m}h_1^{l+m}\ot
   a_{l+m+1}S(h_4^0\cdots h_4^l h_2^{l+1}\cdots h_2^{l+m}))
$$
for all $a_i\in A$ and $h^i\in H$.
This agrees with the cup product on
$\Ext^*_{B^e}(B,B)$ under the given isomorphism.

\medskip

In the special case $A=k$, Theorem \ref{genthm} implies that
\begin{equation}\label{had}
  \HH^*(H)\cong \coh^*(H, H^{ad}),
\end{equation}
where $H^{ad}$ is the $H$-module $H$ under the left adjoint action
defined by $(\ad h)(l) = h_1 l S(h_2)$ for $h,l\in H$,
and $\coh^*(H,H^{ad}):=\Ext^*_H(k,H^{ad})$.
This isomorphism appears in \cite{GK} as Prop.\ 5.6.

\smallskip

\begin{example}{\em
We outline an alternative approach to the Hochschild cohomology of rank
one Hopf algebras that we computed in Section 2, based on
(\ref{had}). This allows us to relate
the structure of $\HH^*(B)$ to the adjoint representation of
$B=H$ for these Hopf algebras. We first find a decomposition of
$B^{ad}$.

Using the coproducts (\ref{coproducts}), we have
$$
  (\ad g)(x^ih) = \chi(g^i)x^ighg^{-1} \ \ \ \mbox{ and } \ \ \
  (\ad x) (x^ih) = (1-\chi(g_1^ih))x^{i+1}h
$$
for all $g,h\in G$. Let $g_0,\ldots,g_t$ be a set of representatives
of conjugacy classes of $G$, where $g_1$ is the central element from
Section 2. Assume $g_1$ has order $m$
 and $g_0=1, \ g_2=g_1^2, \ \ldots , \ g_{m-1}=g_1^{m-1}$.
Then $B^{ad}$ has the following decomposition as a $B$-module
(however the summands are not necessarily indecomposable):
$$
  B^{ad} = \bigoplus_{k=0}^t \Span_k\{x^ihg_kh^{-1}\mid 0\leq i\leq
        n-1, \ h\in G\}.
$$
Each summand above potentially splits into the sum of two
$B$-submodules: For each $h\in G$, let $j_h$ ($0\leq j_h\leq n-1$)
be the smallest such that $\chi(h)=\chi(g_1)^{-j_h}$ if this exists,
and otherwise let $j_h=n-1$. Then the $k$th summand above becomes
$$
 \Span_k\{x^ihg_kh^{-1}\mid 0\leq i\leq j_{g_k}, \ h\in G\}
\oplus \Span_k\{x^ihg_kh^{-1}\mid j_{g_k}+1\leq i\leq n-1, \ h\in
G\},
$$
although this is not needed for the computation of cohomology. In
order to compute $\coh^*(B,B^{ad})$, we use the
free $A$-resolution of $k$:
\begin{equation*}
\cdots \stackrel{\cdot x^{n-1} }{\longrightarrow} A
       \stackrel{\cdot x}{\longrightarrow} A
       \stackrel{\cdot x^{n-1}}{\longrightarrow} A
       \stackrel{\cdot x}{\longrightarrow} A
       \stackrel{\varepsilon}{\rightarrow} k \rightarrow 0.
\end{equation*}
This may be extended to a projective $B$-resolution of $k$ by giving
$A$ the following actions of $G$: In degree $2i$, $\ g\cdot
a=\chi(g)^{in}( {}^{g} a)$, and in degree $2i+1$, $\ g\cdot
a=\chi(g)^{in+1}( {}^{g} a)$ for all $a\in A$ and  $g\in G$. This leads to an
alternative proof of Theorem \ref{gvs}. }\end{example}

\smallskip

We now return to a more general setting.  Assume that $H$ is any
finite dimensional semisimple Hopf algebra and $B=A\# H$.
In this case, a direct consequence of Theorem \ref{genthm} is that
\begin{equation}\label{ss}
\HH^*(B)\cong \HH^*(A,B)^H
\end{equation}
(cf.\ \cite[Thm.\ 3.3]{St}). To see this, first use the relation
$\Hom_{\mathcal D}(M,N)\cong\Hom_{A^e}(M,N)^H$ for any two $\mathcal
D$-modules $M,N$, where the superscript $H$ denotes invariants under
the action
$$
  (h\cdot f)(m) = h_1\cdot (f(S(h_2)\cdot m) = h_1 f(S(h_4)mS^2(h_3))S(h_2)
$$
induced by the embedding $\delta: H\rightarrow H^e$.
Next we
must see that taking $H$-invariants after taking cohomology is
equivalent to taking $H$-invariants before taking cohomology. This
follows from the observation that $H$-invariants are precisely the
image of a nonzero integral since $H$ is semisimple. Using this, we
now give explicit formulas for cocycles and cup products on
$\HH^*(A,B)^H$.

Let
\begin{equation}\label{any}
 \cdots\rightarrow P_2\rightarrow P_1\rightarrow P_0\rightarrow A\rightarrow 0
\end{equation}
be any $\mathcal D$-projective resolution of $A$.
We claim that the bar complex for $A$ is itself a $\mathcal D$-projective
resolution of $A$ where $h\in H\cong \delta(H)$ acts on $A^{\ot m}$ diagonally.
We thank A.\ Kaygun for explaining to us a proof of this fact, in the context
of Hopf-Hochschild cohomology.
We summarize the proof here:
Let $\Lambda$ be an integral for $H$ with $\varepsilon(\Lambda)=1$.
Then the $\mathcal D$-map ${\mathcal D}\rightarrow A^e$ defined by
$ah_1\ot S(h_2)b\mapsto \varepsilon(h)a\ot b$ is split by the $\mathcal D$-map
$A^e\rightarrow {\mathcal D}$ defined by $a\ot b\mapsto a\Lambda_1\ot S(\Lambda_2)b$
for all $a,b\in A$ and $h\in H$.
Therefore $A^e$ is $\mathcal D$-projective.
If $m>0$, we see that $A^{\ot (m+2)}$ is $\mathcal D$-projective as follows.
Since $H$ is semisimple, $A$ is $H$-projective, and so $A^{\ot m}$ is a direct
summand of a sum of copies of $H^{\ot m}$. A standard argument (see the proof
of \cite[Prop.\ 3.1.5]{benson91a}) shows that $H^{\ot m}\cong H\ot (H_{tr}^{\ot (m-1)})$,
where $H_{tr}$ is $H$ with the trivial action of $H$.
Now $A\ot H\ot A\cong {\mathcal D}$ as $\mathcal D$-modules,
and $k$-bases of the remaining factors $H_{tr}^{\ot (m-1)}$ arising from
$A^{\ot (m+2)}$ provide a $\mathcal D$-basis of a free $\mathcal D$-module having
$A^{\ot (m+2)}$ as a direct summand.
It follows that the bar resolution of $A$ is a 
$\mathcal D$-projective resolution of $A$ in case $H$ is semisimple.

Let $\psi_n:A^{\ot (m+2)}\rightarrow
P_m$ be $\mathcal D$-homomorphisms giving a map of chain complexes from the
bar complex  to (\ref{any}). The following theorem
generalizes \cite[Thm.\ 5.4]{CGW}, which is useful in case a
resolution other than a bar-type resolution is used to compute the
cohomology. For example, it was used to find explicit formulas for
Hochschild 2-cocycles in \cite{CGW} when the cohomology was computed
via a Koszul resolution.

\begin{thm}\label{pm}
Assume $H$ is a finite dimensional semisimple Hopf algebra. Let $f:P_m\rightarrow B$ be
a function representing an element of $\HH^m(A,B)^H$ expressed via
the complex (\ref{any}). The corresponding function
$\widetilde{f}\in \Hom_k(B^{\ot m},B)\cong \Hom_{B^e}(B^{\ot
(m+2)},B)$ expressed via the bar complex is defined by
$$
  \widetilde{f}(a_1h^{1}\ot\cdots\ot a_mh^{m}) \hspace{4.3in}
$$
$$\hspace{.2in}=((f\circ\psi_m)(1\!\ot\! a_1\!\ot\! {}^{h_1^{1}} a_2\!\ot\! {}^{h_2^{1}h_1^{2}}
  a_3\!\ot\!\cdots\!\ot {}^{h_{m-1}^{1}h_{m-2}^{2} h_1^{m-1}}a_m\ot 1))
    h_m^{1}h_{m-1}^{2}\cdots h_2^{m-1}h^{m}
$$
for all $a_1,\ldots,a_m\in A$ and $h^{1},\ldots,h^{m}\in H$.
\end{thm}

\begin{proof}
This follows by explicitly tracing through the Eckmann-Shapiro Lemma
as it applies to $\HH^*(B)\cong \Ext^*_{\mathcal D}(A,B)$ in the
proof of Theorem \ref{genthm}. We  use the explicit map from the
bar resolution for $B$ to ${\mathcal D}_* \uparrow_{\mathcal
D}^{B^e}$ given in (\ref{inverse}). We also need a $\mathcal D$-map from
${\mathcal D}_*$ to the bar resolution for $A$, and this is
\begin{equation}\label{further}
  a_0h_1^{0}\!\ot\!\cdots\!\ot\! a_mh_1^{m}\!\ot\! a_{m+1}S(h_2^{0}\cdots
   h_2^{m})\hspace{2.3in}
\end{equation}
$$\hspace{1in}\mapsto
   a_0\!\ot\! {}^{h_1^{0}} a_1\!\ot {}^{h_2^{0}h_1^{1}} a_2\!\ot\!
  \cdots\!\ot\! {}^{h_{m+1}^0 h_m^{1}\cdots h_2^{m-1}h^{m}}a_{m+1}.
$$
(This generalizes \cite[(5.2)]{CGW}.)

Applying (\ref{inverse}) first, $\widetilde{f}(a_1h^{1}\ot\cdots\ot
a_mh^{m}) =\widetilde{f}(1\ot a_1h^{1}\ot\cdots\ot a_mh^{m}\ot 1)$
may be identified with
$$\widetilde{f}((1\!\ot\! h_3^{1}\!\cdots\! h_3^{m})\!\ot\! (1\!\ot\! a_1h_1^{1}
  \!\ot\!\cdots\!\ot\! a_m h_1^{m}\!\ot\! S(h_2^{1}\!\cdots\! h_2^{m}))
\hspace{1in}$$
$$\hspace{1in}=\widetilde{f}(1\!\ot\! a_1h_1^{1}\!\ot\!\cdots\!\ot\! a_mh_1^{m}\!\ot\! S(h_2^{1}\!\cdots
   \! h_2^{m})) h_3^{1}\!\cdots\! h_3^{m}.
$$
Now applying (\ref{further}), this is
$$ ((f\circ\psi_m)(1\!\ot\! a_1\!\ot\! {}^{h_1^{1}} a_2\!\ot\! {}^{h_2^{1}h_1^{2}} a_3
\!\ot\!\cdots\!\ot\! {}^{h_{m-1}^{1}h_{m-2}^{2}\cdots h_1^{m-1}}
a_m\!\ot\! 1))
   h_m^{1}h_{m-1}^{2}\cdots h_2^{m-1}h^{m}.
$$
\end{proof}

We now describe the cup product on $\HH^*(A,B)^H$.
Just as in Section 3, where $H=kG$, the cup product on
$\HH^*(A,B)^H\subseteq \HH^*(A,B)$ induced by the algebra
structure on $B$ corresponds to the cup product on $\HH^*(B)$, and
in particular if $f: A^{\ot (l+2)}\rightarrow B$ and $f':A^{\ot
(m+2)}\rightarrow B$, then
$$
  (f\smile f')(a_0\ot\cdots\ot a_{l+m+1})=f(a_0\ot\cdots\ot a_l\ot 1)
  f'(1\ot a_{l+1}\ot\cdots\ot a_{l+m+1})
$$
by (\ref{further}) and (\ref{cup}).

\section{Hopf-Hochschild cohomology is Hochschild cohomology}

Let $H$ be a bialgebra and $A$ an $H$-module algebra. In \cite{AK},
Kaygun introduces an algebra $\Gamma=A^e \ot H$ with the following
multiplication:
$$(a\ot b\ot h)(c\ot
d\ot l)=(a({}^{h_1}c)\ot ({}^{h_3}d)b\ot h_2l)$$ for all
$a,b,c,d\in A$ and $h,l\in H$. 
By \cite[Lem.\ 3.2]{AK}, $\Gamma$ is an associative
algebra, denoted $A^e\rtimes H$ there. The bar resolution
for $A$ is a differential graded $\Gamma$-module
under the usual action of $A^e$ and the tensor product action of $H$ (see
\cite[Lem.\ 3.5]{AK}).

Let $M$ be an $H$-{\em equivariant} $A$-bimodule, that is $M$ is
both an $H$-module and  an $A$-bimodule, and
$
 {}^{h}(amb) = ({}^{h_1}a)({}^{h_2}m) ({}^{h_3}b)
$
for all $m\in M$, $a,b\in A$, and $h\in H$. Equivalently,
$M$ is a $\Gamma$-module where $(a\ot b\ot h)m =
a({}^{h}m)b$ by \cite[Lem.\ 3.3]{AK}. The {\em
Hopf-Hochschild cohomology} $\HH^*_{\Hopf}(A,M)$ of $A$ with coefficients in $M$ is
defined in \cite{AK} to be the cohomology of the cochain complex
${\rm CH}^*_{\Hopf }(A,M)$ where
$${\rm CH}^m_{\Hopf }(A,M):=\Hom_{\Gamma}(A^{\ot (m+2)},M), $$
and the differentials are induced from those of the bar complex
for $A$.

\begin{lem}\label{gamma}
Let $H$ be a Hopf algebra, and $A$ an
$H$-module algebra. There is an
isomorphism of algebras $\Gamma\cong {\mathcal D}$, where $\mathcal
D$ is defined in (\ref{mathcald}).
\end{lem}

\begin{proof}
Define $\phi: {\mtc D} \ra\Gamma$ by
$$ah_1\ot S(h_2)b \ra a\ot b\ot h$$
for all $a,b\in A$, $h\in H$.
Then clearly $\phi$ has inverse $\psi$ defined by $\psi (a\ot b \ot
h)= ah_1\ot S(h_2)b$. We verify that $\phi$ is multiplicative:
\begin{eqnarray*}
\phi((ah_1\ot S(h_2)b)(cl_1\ot S(l_2)d))&=&
    \phi(a({}^{h_1}c)h_2l_1\ot S(h_3l_2)({}^{h_4}d)b)\\&=&
a({}^{h_1}c)\ot ({}^{h_3}d)b\ot h_2l\\&=& \phi(ah_1\ot
S(h_2)b)\phi(cl_1\ot S(l_2)d)
\end{eqnarray*}
for all $a,b,c,d\in A$ and $h,l\in H$.
\end{proof}

If $M$ is an $A\# H$-bimodule, then $M$ is a module for $\Gamma\cong {\mathcal D}$
by restriction to ${\mathcal D}\subset (A\# H)^e$.
Therefore $M$ has the structure of an
$H$-equivariant $A$-bimodule.
The next theorem shows that the Hopf-Hochschild cohomology of $A$ with
coefficients in $M$ is isomorphic to Hochschild cohomology under the
assumption that $H$ is semisimple (cf.\ \cite[Thm.\ 3.7]{AK}).

\begin{thm}\label{Hopf-Hochschild}
Let $H$ be a finite dimensional semisimple Hopf algebra, and $A$ an
$H$-module algebra. Then 
$$
  \HH^m_{\Hopf}(A,M)\cong \HH^m(A\# H, M)
$$
for all $m$ and any $A\# H$-bimodule $M$.
\end{thm}

\begin{proof}
The bar resolution of $A$ is a $\mathcal D$-projective resolution
of $A$, as explained towards the end of Section 4. This fact, together with
Lemma \ref{gamma} and Theorem \ref{genthm} imply
$$\HH^m_{\Hopf}(A,M)\cong\Ext_{\mtc D}^m(A, M)\cong\HH^m(A\#H, M).$$
\end{proof}

\end{document}